\documentclass[]{article}

\usepackage[T1]{fontenc}
\usepackage[frenchb]{babel}
\usepackage{amsmath}
\usepackage{amssymb}
\usepackage{enumerate}
\begin{document}


\title{Henri Poincar\'e and his ''model'' of hyperbolic geometry}

\author{Philippe Nabonnand\\
Laboratoire d'Histoire des Sciences et de Philosophie \\CNRS -- Universit\'e de Lorraine}
\maketitle
The aim of the talk is to trace how and when Henri Poincar\'e used non-Euclidean geometries (NEG) in his mathematical and philosophical works, with a particular attention to the genesis and the description of his ''model''. We begin by a short presentation of the context of NEG in France around the  1870-80s. Then we expound from several sources the introduction and use of NEG in  Poincar\'e's work about Fuchsian functions and we stress on the analogy between elliptic functions and fuchsian functions.
\section{The context of non-Euclidean geometry in France around the years 1870-80s}
At the end of 1869, Jules Carton sent to the Academy of Sciences of Paris a ''proof'' of the postulate of parallels. One of the leaders of the Academy, the well-known mathematician Joseph Bertrand, approved this proof. This announcement and this almost official approval provoked a series of proposals for proof of the postulatum but also many criticisms, first, expressed in the privacy of correspondence, but then quickly in newspapers such as the review published by the Abbot Moigno, \emph{Cosmos}.

Others, like Jules Ho\"uel and Gaston Darboux, saw it as an opportunity to popularize and deepen the debate. Jules Ho\"uel was the translator in French of the major texts of non-Euclidean geometries. Ho\"uel fought for the acceptance and the recognition of NEG in a context of discussions about the provability of the axiom of parallels, the consistency of NEG and the status of the axioms of geometry. His point of view was moderately empiricist.

Since 1875, there had been a reception of NEG and a debate in the field of Philosophy via the \emph{Revue philosophique de la France et de l'\'Etranger}. Founded by the psychologist Th\'eodule Ribot, the \emph{Revue} gave special attention to contemporary debates on philosophy of science, with a focus on NEG and the status of axioms of geometry. Many actors contributed to this debate, mathematicians, engineers, psychologists-physiologists. In this context, the \emph{Revue} stressed the importance of German theories in experimental psychology, especially about ''spatial sense''.

\section{The three ''Suppl\'ements''}
In 1880, Henri Poincar\'e took part in a competition announced by the French Academy in 1878. The subject was
''To perfect in any material respect the linear differential equations theory with a single independent variable'.
First he submitted a memoir to which he later added three ''suppl\'ements\footnote{These three ''suppl\'ements'' were discovered in the Archives of Academy of Sciences of Paris by J.~J.~Gray in 1995 and edited, with an introduction, in 1997 by J.~J.~Gray and S.~Walter \cite{GW97}.}''.

In the first one, Poincar\'e studied the behavior of the quotient $z=\frac{f(x)}{g(x}$ of two independent solutions of a linear differential equation of order 2 and asked the question to know when $x$ is a meromorphic function of  $z$. In this intention, he described a subgroup of transformations of $PGL(2,\mathbb{R})$ and an associated tessellation (paving) of the unity disk. Poincar\'e stressed the link of these geometrical considerations with the hyperbolic Geometry.  For this goal, he identified the group of transformations he studied and the group of the ''pseudogeometry'' of Lobatchewski. In fact, he will made a very moderate use of the ''convenient language'' but at the end of the first supplement he introduced a seminal remark which he would thereafter consider as the core of the use of NEG in the theory of Fuchsian functions, the analogy between elliptic functions and Fuchsian functions.

In the \emph{Report on his own works} [\cite{HP86}, he explains the crucial nature of the use of NEG in the theory of Fuchsian functions as resulting from the analogy elliptic functions/Fuchsian functions. The analogy breaks down as follows:
\begin{center}
\begin{tabular}{cccc}
Euclidean  & Discrete subgroups 	& Lattices &	Elliptic \\
geometry & of orthogonal group & &functions\\[2ex]
Non-Euclidean 	& Discrete subgroups  &	Hyperbolic 	& Fuchsian\\ 
geometry& of $PSL(2,\mathbb{R})$ & pavings & functions
\end{tabular}
\end{center}

\noindent\emph{Abstract}:
	\begin{enumerate}[.~]
	\item Few drawings considered in the context of geometry of the unity disk.
	\item Identification of the groups of transformation $=$ identification of the geometries.
	\item The identification of geometries provides a convenient language.
	\item The thema of the analogy elliptic functions/Fuchsian functions.
	\end{enumerate}

\noindent In the second supplement, Poincar\'e gave a definition of the elements of pseudo geometric plane in terms of classical geometry of disc unity. He described also the group of pseudo geometric movements in terms of homographies which set the fundamental circle.

\section{The first half of 1881}
In a talk about ''applications of NEG to theory of quadratic forms'' \cite{HP81a},
 Poincar\'e uses, despite the title, the same exposition mode as in the ''suppl\'ements''. He first studies the linear transformations (with integral coefficients) which preserve a ternary quadratic form (with integral coefficients). Following Hermite and Selling, he is led to investigate the geometry of tessellations of unity disc. After a classical description of the geometry of the group of substitutions that exchange regions of the tessellation, he finds it convenient to use the vocabulary of the pseudo geometry.

\noindent\emph{Abstract}:
\begin{enumerate}[.~]
	\item Identification of geometries as identification of elements  .
	\item The identification of geometries provides a convenient language.
\end{enumerate}

Poincar\'e published eight notes about Fuchsian functions  during the first half of the year; only three of them mention NEG. This raises a question: are NEG really important for Poincar\'e's theory of Fuchsian functions? Poincar\'e's answer is ambivalent. He emphizes that NEG are very important for the discovery process but he doesn't really use NEG in his papers \cite{HP81b}.

In a note about Kleinian groups  (published July 11th, 1881) \cite{HP81c}, Poincar\'e copes with the question of finding discrete subgroups of $PSL(2,\mathbb{C})$. Of course, finding Kleinian groups is a more general problem than finding Fuchsian groups, which are discrete subgroups of $PSL (2,\mathbb{R})$. Once again, Poincar\'e explains how NEG are important in the discovery process without translating it explicitly in the exposition of theory. In this paper, he gives a description of hyperbolic geometry on a half-space (3-dimensional hyperbolic geometry).
In his paper on Fuchsian groups in \emph{Acta mathematica}, Poincar\'e evokes NEG in the same terms.

\emph{Abstract}:
\begin{enumerate}[.~]
\item A claim that NEG was important for the discovery of Fuchsian and Kleinian groups.
\item A new identification of elements.
\item No real use of NEG
\end{enumerate}
\section{Conclusion}
In a paper entitled 'Les g\'eom\'etries non euclidiennes' \cite{HP91}, Poincar\'e claims that his dictionary is a proof of the non-contradiction of hyperbolic geometry.
In this context, we can say that the half plane of Poincar\'e is a model (in the logical sense\footnote{If a deductive system has a model, the system is semantically consistent.}) but we have to notice that the translation of axioms of NEG is not explicite (perhaps, included in the claim concerning all the theorems).

In any case, Poincar\'e made a very moderate use or no-use of the ''convenient language'' in mathematical papers. In particular, there is no drawing when dealing with NEG. Nevertheless, refering to the analogy between elliptic functions and Fuchsian functions, he claimed that hyperbolic geometry played a crucial role in the process of discovery.

Following the differentiation between structural analogy (correspondence between relations) and functional analogy  (correspondence between elements which have analogous properties), we can notice that the functional part of the correspondence in Poincar\'e's dictionnary of Poincar\'e is explicit and that the functional part is implicit (excepted when Poincar\'e refers to isomorphism between groups); nevertheless, Poincar\'e's conclusions (correspondence between theorems) are true if the analogy is structural.


\end{document}